\documentclass{elsart3-1_intemporel}
\usepackage{amssymb}
\usepackage{xspace}
\usepackage{mathrsfs}
\usepackage[latin1]{inputenc}  % accents 8 bits dans la source
\usepackage[T1]{fontenc}       % accents dans le DVI
\usepackage[english,francais]{babel}

\newtheorem{theorem}{Theorem}[section]

\newtheorem{e-proposition}[theorem]{Proposition}
\newtheorem{corollary}[theorem]{Corollary}
\newtheorem{e-definition}[theorem]{Definition\rm}

\renewcommand{\theequation}{\arabic{equation}}
\newcommand{\RR}{\mathbf{R}}

\newcommand{\ZZ}{\mathbf{Z}}

\newcommand{\SL}{\mathrm{SL}}

\newcommand{\norma}{\mathscr{N}}
\newcommand{\centra}{\mathscr{Z}}
\newcommand{\QZ}{\mathscr{Q\!Z}}
\newcommand{\se}{\subseteq}

\newcommand{\cat}{{\upshape CAT(0)}\xspace}
\newcommand{\tangle}[2]% angle de Tits
{\angle_\mathrm{T}(#1,#2)}
\newcommand{\aangle}[3]% angle d'Alexandrov
{\angle_{#1}(#2,#3)}
\newcommand{\cangle}[3]% angle de comparaison
{\overline{\angle}_{#1}(#2,#3)}

\newcommand{\Isom}{\mathrm{Is}}
\newcommand{\soc}{\mathrm{soc}}
\newcommand{\bd}{\partial}

\def\og{\leavevmode\raise.3ex\hbox{$\scriptscriptstyle\langle\!\langle$~}}
\def\fg{\leavevmode\raise.3ex\hbox{~$\!\scriptscriptstyle\,\rangle\!\rangle$}}

\journal{the Acad\'emie des sciences}
\begin{document}
% place in the next line the header (rubrique) chosen for your article,
% if you know it (you can also have 2, format : Header1/Header2
\centerline{}
\begin{frontmatter}


% \author{Name\thanksref{label2}}
\selectlanguage{english}
\title{Some properties of non-positively curved lattices}
% The [label1] can be suppressed if there is only one address for all authors

\selectlanguage{english}
\author[PEC]{Pierre-Emmanuel Caprace\thanksref{thankPEC}},
\ead{caprace@ihes.fr}
\author[NM]{Nicolas Monod\thanksref{thankNM}}
\ead{nicolas.monod@epfl.ch}

\address[PEC]{IHES, France}
\address[NM]{EPFL, Switzerland}
\thanks[thankPEC]{Supported in part by IPDE and a Hodge fellowship.}
\thanks[thankNM]{Supported in part by FNS.}
% If you know the dates of reception, and acceptation you can put them now;
%  idem the name of the person presenting the Note

\medskip %                      ATTENTION, suppression :
%\begin{center}
%{\small Received *****; accepted after revision +++++\\
%Presented by £££££}
%\end{center}

\begin{abstract}
\selectlanguage{english} We announce results on the structure of \cat groups, \cat lattices and of the
underlying spaces. Our statements rely notably on a general study of the full isometry groups of proper \cat
spaces. Classical statements about Hadamard manifolds are established for singular spaces; new arithmeticity
and rigidity statements are obtained.
%{\it To cite this article: P.-E. Caprace, N. Monod, C. R. Acad. Sci. Paris, Ser. I 340 (2005).}

\vskip 0.5\baselineskip

\selectlanguage{francais} \noindent{\bf R\'esum\'e} \vskip 0.5\baselineskip \noindent {\bf Quelques
propri\'et\'es des groupes \cat. } Nous pr\'esentons des r\'esultats de structure sur les groupes \cat, les
r\'eseaux \cat et sur les espaces sous-jacents. Nos \'enonc\'es reposent notamment sur une \'etude g\'en\'erale
des groupes d'isom\'etries pleins des espaces \cat propres. Nous démontrons des résultats qui généralisent des
énoncés classiques sur les variétés de Hadamard et proposons de nouveaux théorèmes d'arithméticité et rigidité.
%{\it Pour citer cet article: P.-E. Caprace, N. Monod, C. R. Acad. Sci. Paris, Ser. I 340 (2005).}
\end{abstract}
\end{frontmatter}

% now the Version française abr\'eg\'ee, if it exists
\selectlanguage{francais}
\section*{Version fran\c{c}aise abr\'eg\'ee}
% Text of your Version française abr\'eg\'ee here.
% Note you do not need to repeat here equations that you use in the
% main text - for example 'voir (3)' is quite acceptable.

Nous conviendrons qu'un \textbf{groupe \cat} est un couple $(\Gamma, X)$ où $X$ est un espace \cat propre et
$\Gamma$ un groupe d'isométries de $X$ dont l'action est proprement discontinue et cocompacte. Il s'agit
d'étudier les relations entre la géométrie de $X$ et les propriétés algébriques de $\Gamma$. Ce cadre permet un
traitement unifié de nombreuses situations classiques (groupes fondamentaux de variétés compactes à courbure
négative, réseaux uniformes des groupes algébriques semisimples, en particulier groupes $S$-arithmétiques
anisotropes) et moins classiques liées à la théorie géométrique des groupes (réseaux non linéaires associés aux
arbres, immeubles exotiques et non euclidiens, nombreux groupes Gromov-hyperboliques).
%
%\smallskip
Cette note annonce quelques résultats généraux qui relèvent de ce contexte. Nous considérons parfois le cas plus
général des \textbf{réseaux \cat} qui, à défaut de meilleure définition, consistent des couples $(\Gamma, X)$
formés d'un réseau $\Gamma$ dans un groupe d'isométries cocompact de $X$.

\smallskip Voici un premier exemple des relations entre $X$ et $\Gamma$: la dimension du facteur euclidien de
$X$ est égale au rang maximal des sous-groupes abéliens libres normaux dans $\Gamma$
(\textbf{Cor.~\ref{cor:lattices:EuclideanSplitting}}; ces numéros renvoient à la partie anglaise). Dans le cas
des variétés riemanniennes, c'est là un résultat d'Eberlein~\cite{Eberlein83} qui peut être vu comme une
réciproque (partielle) au \og théorème du tore plat\fg~\cite{Gromoll-Wolf}, \cite{Lawson-Yau},
\cite[\S II.7]{Bridson-Haefliger}.

\smallskip Pour obtenir des énoncés de rigidité, il est utile (et  bien souvent nécessaire) de supposer que $X$
soit \textbf{géodésiquement complet}, ce qui n'exclut aucun des exemples classiques (immeubles, variétés,
\ldots). Nous montrons alors que $X$ possède une isométrie parabolique si et seulement si $X$ admet une
décomposition isométrique $X=M\times X'$ où $X'$ est \cat et $M$ est un espace symétrique de type non compact.
(\textbf{Cor.~~\ref{cor:parabolic}}; cet énoncé est faux pour des espaces $X$ sans réseau $\Gamma$). Plus
précisément, soit $(\Gamma, X)$ un groupe \cat avec $\Gamma$ irréductible et $X$ géodésiquement complet.
Supposons que $X$ admette une isométrie parabolique. Si $\Gamma$ est résiduellement fini, alors $X$ est
isométrique à un produit d'espaces symétriques et d'immeubles de Bruhat--Tits. Sinon, l'intersection $\Gamma_D $
 de tous les sous-groupes d'indice fini de $\Gamma$ n'est pas de type fini et le quotient $\Gamma/\Gamma_D$
 est un groupe arithmétique
(\textbf{Thm.~\ref{thm:arith:geometric}}).

\smallskip Notre travail s'appuie sur une analyse des groupes d'isométries des espaces \cat propres indépendante
de l'existence de réseaux. On se place dans le cadre où un groupe $G<\Isom(X)$ agit \textbf{minimalement} et
sans point fixe à l'infini (il convient de montrer qu'il est possible de se restreindre à ce cas). Lorsque le
bord à l'infini de $X$, muni de la métrique de Tits, est de dimension finie, on montre que $X$ possède une
décomposition canonique en un produit d'un facteur euclidien et d'un nombre fini de facteurs irréductibles non
euclidiens; c'est là une variante de la décomposition de de Rham obtenue dans \cite{FoertschLytchak06}. En
particulier le groupe d'isométries complet de $X$ se décompose virtuellement comme produit des groupes
d'isométries de chaque facteur de $X$. En outre, le groupe d'isométries de tout facteur irréductible non
euclidien est soit un groupe de Lie simple presque connexe de centre trivial, soit un groupe totalement
discontinu dont le radical moyennable est trivial. Dans tous les cas, ce groupe est irréductible en ce sens
qu'aucun sous-groupe fermé d'indice fini ne se scinde en un produit direct non trivial
(\textbf{Thm.~\ref{thm:decomposition}}). En fait, tout groupe d'isométries d'un espace \cat irréductible non
euclidien  dont le bord est de dimension finie transmet à chacun de ses sous-groupes normaux non triviaux la
propriété d'agir minimalement et sans point fixe au bord sur l'espace en question
(\textbf{Thm.~\ref{thm:irred}}). De cette propriété de \og densité géométrique \fg des sous-groupes normaux, que
l'on peut interpréter comme une forme faible de simplicité, découlent des énoncés purement algébriques: un
sous-groupe normal, et plus généralement un sous-groupe sous-normal, ne se scinde pas en produit direct non
trivial, son radical moyennable et son centralisateur sont tous deux triviaux. Lorsque non seulement le groupe
d'isométries $\Isom(X)$, mais aussi chacun de ses sous-groupes ouverts, agit sans fixer de point à l'infini, ces
différentes conclusions peuvent être considérablement renforcées (\textbf{Thm.~\ref{thm:NoOpenStabiliser}}).

\smallskip
Le phénomène de \og densité géométrique \fg qu'on vient de décrire pour les sous-groupes normaux d'un groupe
d'isométries d'un espace \cat irrédutible se manifeste également pour les réseaux
(\textbf{Thm.~\ref{thm:density}}). Le terme de \og densité \fg est ici particulièrement approprié: de ces
considérations, on déduit en outre un preuve du résultat classique, dû à Borel, de densité de Zariski dans le
cas particulier des réseaux de groupes semisimples.

\smallskip
Mentionnons finalement que ces différents développements permettent d'appliquer le théorème de superrigidité de
Margulis pour certains groupes arithmétiques tels que $\SL_n(\ZZ)$ dans un contexte purement \cat
(\textbf{Thm.~\ref{thm:superrigidity}}). Par ailleurs, on montre qu'au sein des espaces \cat propres
géodésiquement complets, les espaces symétriques et immeubles euclidiens fortement homogènes sont caractérisés
par le fait que le fixateur de tout point à l'infini est cocompact
(\textbf{Thm.~\ref{thm:cocompact:stabilisers}}).

\newpage
\selectlanguage{english}
% main text
%\section{Introduction}
%%%%%%%%%%%%%%%%%%%%%%%%%%%%%%%%%%%%%%%%%%%%%%%%%%%%%%%%%%%%%%%%%%%%%%%%%%%%%%%%%%%%%%%%%%%%%%%%%%%%%%%%%%%%%%
\section{\cat groups and lattices}

We define a \textbf{\cat group} as a pair $(\Gamma, X)$ where $X$ is a proper \cat space and $\Gamma$ a properly
discontinuous cocompact group of isometries of $X$. The study of these objects centers around the following
general question:

\smallskip
\begin{center}
\itshape
$(\star)$ What is the interplay between the geometry of $X$ and the algebraic properties of $\Gamma$?
\end{center}
\smallskip

The motivation for considering \cat groups is that they provide a common framework for at least four classical
topics: closed manifolds of non-positive curvature; uniform lattices in semi-simple groups over local fields, in
particular anisotropic $S$-arithmetic groups; non-linear cousins such as tree lattices, exotic and non-Euclidean
buildings; general geometric group theory. Actually, the classical situations alluded to above naturally suggest
to relax the cocompactness condition in order to cover \emph{all} lattices in semi-simple groups, as well as
larger families of non-linear relatives, including non-uniform lattices arising from Kac--Moody theory. In the
general \cat setting, there is (as yet) no consensus on a good definition for discrete groups of finite
covolume; we shall content ourselves with the following \emph{ad hoc} definition:

A \textbf{\cat lattice} is a pair $(\Gamma, X)$ where $X$ is a proper \cat space whose isometry group $\Isom(X)$
acts cocompactly and $\Gamma$ is a lattice in $\Isom(X)$. We emphasise that a \cat group is in particular a
finitely generated (uniform) \cat lattice.

\smallskip
The goal of this note is to announce a few general results on \cat lattices, of relevance to each of the above themes.

\medskip
A first example towards question~$(\star)$ regards the maximal Euclidean factor. We recall that the \emph{Flat
Torus theorem}, originating in the work of Gromoll--Wolf~\cite{Gromoll-Wolf} and Lawson--Yau~\cite{Lawson-Yau},
associates Euclidean subspaces $\RR^n\se X$ to subgroups $\ZZ^n \se\Gamma$ (see~\cite[\S II.7]{Bridson-Haefliger}).
The converse is a well-known open problem (Gromov~\cite[\S$6.\mathrm{B}_3$]{Gro93}; for manifolds see Yau,
problem~65 in~\cite{YauPB}). We propose the following.

\smallskip
\begin{theorem}\label{thm:lattices:EuclideanSplitting}
Let $(\Gamma, X)$ be a finitely generated \cat lattice and let $X\cong\RR^n\times X'$ be the Euclidean
decomposition. Then there is a finite index subgroup $\Gamma_0$ which splits as a direct product $\Gamma_0 \cong
\ZZ^n \times \Gamma'$.
\end{theorem}
\smallskip

In the special case of cocompact Riemannian manifolds, the above was the main result of Eberlein's
article~\cite{Eberlein83}. A \cat space $X$ is called \textbf{minimal} if $\Isom(X)$ acts minimally, i.e.
without stabilising any non-empty closed convex proper subset.

\smallskip
\begin{corollary}\label{cor:lattices:EuclideanSplitting}
If $X$ is minimal $($\emph{e.g.} geodesically complete$)$, then the dimension of the Euclidean factor of $X$ equals
the maximal $\mathbf{Q}$-rank of an Abelian normal subgroup of $\Gamma$.
\end{corollary}
\smallskip

In the sequel, an abstract group $\Gamma$ is called \textbf{irreducible} if no finite index subgroup splits as a
non-trivial direct product. We generalise to all finitely generated \cat lattices the Margulis
irreducibility criterion for lattices in semi-simple groups.

%%%%%%%%%%%%%%%%%%%%%%%%%%%%%%%%%%%%%%%%%%%%%%%%%%%%%%%%%%%%%%%%%%%%%%%%%%%%%%%%%%%%%%%%%%%%%%%%%%%%%%%%%%%%%%
\section{Geometric arithmeticity}
A \textbf{neutral} parabolic isometry is a fixed-point free isometry with zero translation length.

\smallskip
\begin{theorem}\label{thm:ArithmeticityCAT0}
Let $(\Gamma, X)$ be an irreducible \cat group. If $X$ admits any neutral parabolic isometry, then either:
\begin{enumerate}
\item $\Isom(X)$ is a rank one simple Lie group with trivial centre; or:

\item $\Gamma$ possesses a normal subgroup $\Gamma_D$ such that $\Gamma/\Gamma_D$ is an arithmetic group. Moreover,
$\Gamma_D$ is either finite or infinitely generated.
\end{enumerate}
\end{theorem}
\smallskip

Assuming in addition that every geodesic segment of $X$ can be extended to a bi-infinite geodesic line (which need not
be unique), we obtain geometric information on $X$ and at the same time drop the neutrality assumption.

\smallskip
\begin{theorem}\label{thm:arith:geometric}
Let $(\Gamma, X)$ be an irreducible \cat group with $X$ geodesically complete. Assume that $X$ possesses some
parabolic isometry.

If $\Gamma$ is residually finite, then $X$ is a product of symmetric spaces and Bruhat--Tits buildings. In
particular, $\Gamma$ is an arithmetic lattice unless $X$ is a real or complex hyperbolic space.

If $\Gamma$ is not residually finite, then $X$ still splits off a symmetric space factor. Moreover, the finite
residual $\Gamma_D$ of $\Gamma$ is infinitely generated and $\Gamma/\Gamma_D$ is an arithmetic group.
\end{theorem}
\smallskip

\begin{corollary}\label{cor:parabolic}
Let $(\Gamma, X)$ be a \cat group with $X$ geodesically complete. Then $X$ possesses a parabolic isometry if and
only if  $X \cong M \times X'$, where $M$ is a symmetric space of non-compact type.
\end{corollary}
\smallskip

For lattices in products of groups that are \emph{simple}, or have few factors, an arithmeticity/non-linearity
alternative was proved in~\cite{Monod_alternative}. In our geometric setting, we can establish it without any
assumption on the factors, and moreover establish geometric superrigidity.

\smallskip
\begin{theorem}\label{thm:arith:geometric:lin}
Let $(\Gamma, X)$ be an irreducible \cat group with $X$ geodesically complete. Assume that $\Gamma$ possesses
some faithful finite-dimensional linear representation $($in characteristic~$\neq 2,3$$)$.

If $X$ is reducible, then $\Gamma$ is an arithmetic lattice and $X$ is a product of symmetric spaces and
Bruhat--Tits buildings.
\end{theorem}

%%%%%%%%%%%%%%%%%%%%%%%%%%%%%%%%%%%%%%%%%%%%%%%%%%%%%%%%%%%%%%%%%%%%%%%%%%%%%%%%%%%%%%%%%%%%%%%%%%%%%%%%%%%%%%
\section{A geometric Borel density theorem}
A recurring theme of our work is that \emph{minimality}~--- a much weaker assumption, upon adjustments, than
the familiar notions of cocompactness or of full limit sets~--- is a valuable geometric notion, similar to Zariski density
in algebraic groups. The following corresponds to Borel's classical result~\cite{Borel60} (and contains it, indeed).

\smallskip
\begin{theorem}\label{thm:density}
Let $G$ be a locally compact group with a continuous isometric  action on a proper \cat space $X$ without
Euclidean factor.

If $G$ acts minimally and without global fixed point in $\bd X$, then any closed subgroup with finite
invariant covolume in $G$ still has these properties.
\end{theorem}
\smallskip

One deduces generalisations of some facts known in the case of lattices in semi-simple groups.

\smallskip
\begin{corollary}
\label{cor:LatticeNormaliser}
Let $X$ be a proper \cat space without Euclidean factor such that $G = \Isom(X)$ acts  minimally without fixed
point at infinity, and let $\Gamma < G$ be a closed subgroup with finite invariant covolume. Then:
\begin{enumerate}
\item $\Gamma$ has trivial amenable radical.\label{pt:LatticeNormaliser:ramen}

\item The centraliser $\centra_G(\Gamma)$ is trivial.\label{pt:LatticeNormaliser:centra}

\item If $\Gamma$ is finitely generated, then is has finite index in its normaliser $\norma_G(\Gamma)$
and the latter is a finitely generated lattice in $G$.
\end{enumerate}
\end{corollary}
\smallskip

As pointed out by P. de la Harpe, (ii)~implies in particular that any lattice in $G$ is~ICC and hence its von Neumann algebra is a factor.

%%%%%%%%%%%%%%%%%%%%%%%%%%%%%%%%%%%%%%%%%%%%%%%%%%%%%%%%%%%%%%%%%%%%%%%%%%%%%%%%%%%%%%%%%%%%%%%%%%%%%%%%%%%%%%
\section{Isometry groups and their normal subgroups}
Our results on \cat groups and lattices require some groundwork on the full isometry group of the underlying \cat spaces.
A common first step for many of our result is the following fact.

\smallskip
\begin{theorem}\label{thm:ramen}
Let $X$ be a proper \cat space with finite-dimensional boundary and no Euclidean factor. Let
$G<\Isom(X)$ be a closed subgroup acting minimally and without fixed point at infinity. Then the amenable radical
of $G$ is trivial.
\end{theorem}
\smallskip

Next, we establish a group decomposition, supplemented by a de Rham decomposition of the space which is a variant of~\cite{FoertschLytchak06}.

\smallskip
\begin{theorem}\label{thm:decomposition}
Let $X$ be a proper minimal \cat space with finite-dimensional boundary. If $G = \Isom(X)$ has no global fixed
point at infinity, then  there is a canonical finite index open characteristic subgroup $G^*\lhd \Isom(X)$ which
admits a canonical decomposition
$$G^*  \cong\ S_1\times \cdots \times S_p \times \big(\RR^n\rtimes
\mathbf{O}(n)\big) \times D_1\times \cdots \times D_q  \kern10mm(p,q,n\geq 0)$$
where $S_i$ are almost connected simple Lie groups with trivial centre and $D_j$ are totally disconnected
irreducible groups. Furthermore, there is a canonical equivariant isometric splitting
$$X \cong\ X_1\times \cdots \times X_p \times \RR^n \times Y_1\times \cdots \times Y_q$$
with componentwise minimal action; all $X_i$ and $Y_j$ are irreducible. Any other product decomposition of $G^*$
or $X$ is a regrouping of the above factors.
\end{theorem}
\smallskip

The conclusion of Theorem~\ref{thm:decomposition} can be considerably strengthened by a further description of
the factors $D_i$. Indeed, they satisfy the following geometric form of simplicity.

\smallskip
\begin{theorem}\label{thm:irred}
Let $X\neq \RR$ be an irreducible proper \cat space with finite-dimensional Tits boundary and $G<\Isom(X)$ any
subgroup whose action is minimal and does not have a global fixed point in $\bd X$.

Then every non-trivial normal subgroup $N\lhd G$ still acts minimally and without fixed point in $\bd X$.
Moreover, the amenable radical of $N$ and its centraliser $\centra_{\Isom(G)}(N)$ are both trivial; $N$ does not
split as a product.

These conclusions hold more generally when $N<G$ is any non-trivial subnormal or even ascending subgroup.
\end{theorem}
\smallskip

Notice that Theorem~\ref{thm:irred} can be gainfully combined with Theorem~\ref{thm:density} and
Corollary~\ref{cor:LatticeNormaliser} in order to obtain for instance information on normal subgroups of
lattices, or lattices in normal subgroups.

\medskip

For the study of the totally disconnected factors of isometry groups, the following \emph{smoothness} result
is a basic link between the topology and the algebra.

\smallskip
\begin{theorem}\label{thm:smooth}
Let $X$ be a geodesically complete proper \cat space and $G<\Isom(X)$ a totally disconnected closed subgroup
acting minimally. Then the pointwise stabiliser in $G$ of every bounded set is open.
\end{theorem}
\noindent
(The conclusion can indeed fail when $X$ is not geodesically complete.)
\smallskip

At this point, we notice in passing that we have gathered enough information about general cocompact \cat spaces
to apply Margulis' superrigidity at least for some lattices.

\smallskip
\begin{theorem}\label{thm:superrigidity}
Let $X$ be a proper \cat space whose isometry group acts cocompactly and without global fixed point at infinity.
Let $\Gamma = \SL_n(\ZZ)$ with $n \geq 3$ and $G = \SL_n(\RR)$.

For any isometric $\Gamma$-action on $X$ there is a non-empty $\Gamma$-invariant closed convex subset $Y\se X$
on which the $\Gamma$-action extends uniquely to a continuous isometric action of $G$.

\smallskip\nobreak\noindent
$($The corresponding statement applies to all those lattices in semi-simple Lie groups that have virtually
bounded generation by unipotents. It also applies to $S$-arithmetic lattices such as
$\SL_n(\ZZ[{\scriptstyle\frac{1}{m}}])$, $n\geq 3$.$)$
\end{theorem}
\smallskip

Observe that the above theorem has no assumptions whatsoever on the action; cocompactness is an assumption on
the given \cat space. It may happen that $\Gamma$ fixes points in $\bd X$, but its action on $Y$ is always
minimal and without fixed points at infinity.

\smallskip

Recall that the \textbf{quasi-centre} $\QZ(G)$ of a group $G$ is the union of the discrete conjugacy classes; it
is always a (topologically) characteristic subgroup. The following result establishes in particular the
existence of a minimal non-trivial normal subgroup; to state it, we define the \textbf{socle} $\soc(\cdot)$ as
the subgroup generated by the (possibly empty) collection of all minimal non-trivial closed normal subgroups.

\smallskip
\begin{theorem}\label{thm:NoOpenStabiliser}
Let $X$ be a proper geodesically complete \cat space without Euclidean factor and $G < \Isom(X)$ a closed subgroup acting cocompactly.
Suppose that no open subgroup of $G$ fixes a point at infinity.% Then we have the following:
\begin{enumerate}
\item $X$ admits a canonical equivariant splitting
$$X \cong\ X_1\times \cdots \times X_p \times  Y_1\times \cdots \times Y_q$$
where each $X_i$ is a symmetric space and each $Y_j$ possesses a $G$-invariant locally finite decomposition
into compact convex cells.

%\item $G$ possesses hyperbolic elements.
%
\item Every compact subgroup of $G$ is contained in a maximal one; the maximal compact subgroups fall into
finitely many conjugacy classes.

\item $\QZ(G) = \{1\}$; in particular $G$ has no non-trivial discrete normal subgroup.

\item $\soc(G^*)$ is a direct product of $p+q$ non-discrete characteristically simple groups.
\end{enumerate}
\end{theorem}
\smallskip

Finally, we present a result of a different vein. Symmetric spaces and Bruhat--Tits buildings have in common the
property that the stabilisers of points at infinity are cocompact (being always parabolic; the case of Euclidean
type is obvious). This property is further shared by \emph{Bass--Serre} trees, \emph{i.e.}  edge-transitive
trees (which are in particular regular or bi-regular).

\smallskip
\begin{theorem}\label{thm:cocompact:stabilisers}
Let $X$ be a geodesically complete proper \cat space. Suppose that the stabiliser of every point at infinity
acts cocompactly on $X$.

Then $X$ is isometric to a product of symmetric spaces, Euclidean buildings and Bass--Serre trees.
\end{theorem}
\smallskip

The Euclidean buildings appearing in the preceding statement admit an automorphism group that is strongly
transitive, \emph{i.e.} acts transitively on pairs $(c, A)$ where $c$ is a chamber and $A$ an apartement
containing $c$. This property characterises the Bruhat--Tits buildings, except perhaps for some two-dimensional
cases where this is a known open question.

\end{document}